\documentclass[conf]{new-aiaa}
\usepackage[utf8]{inputenc}

\usepackage{graphicx}
\usepackage{amsmath}
\usepackage[version=4]{mhchem}
\usepackage{siunitx}
\usepackage{longtable,tabularx}
\usepackage{multirow}
\usepackage{subcaption}
\usepackage{cleveref}
\usepackage{algorithm}
\usepackage{algpseudocode}

\newcolumntype{M}[1]{>{\centering\arraybackslash}m{#1}}
\usepackage{tikz,pgfplots}
\pgfplotsset{compat=1.12}
\usepackage{pgfplotstable}
\usepackage{tablefootnote}
\pgfplotsset{
    discard if/.style 2 args={
        x filter/.append code={
            \edef\tempa{\thisrow{#1}}
            \edef\tempb{#2}
            \ifx\tempa\tempb
                
            \fi
        }
    },
    discard if not/.style 2 args={
        x filter/.append code={
            \edef\tempa{\thisrow{#1}}
            \edef\tempb{#2}
            \ifx\tempa\tempb
            \else
                
            \fi
        }
    }
}
\definecolor{hcorange}{RGB}{245, 130, 48}
\definecolor{hcnavy}{RGB}{0, 0, 128}
\definecolor{hcblue}{RGB}{0, 130, 200}
\definecolor{hcpink}{RGB}{250, 190, 190}
\definecolor{hcyellow}{RGB}{255, 255, 255}
\definecolor{hcbrown}{RGB}{128, 0, 0}
\definecolor{hclavender}{RGB}{230, 190, 255}
\definecolor{hcgrey}{RGB}{128, 128, 128}
\definecolor{hcgreen}{RGB}{60, 180, 75}
\definecolor{hcred}{RGB}{230, 25, 75}
\definecolor{green1}{cmyk}{1,0,1,0}

\title{Towards a Scalable Hierarchical High-order CFD Solver}

\author{Zan Xu\footnote{Ph.D. Candidate, Department of Aeronautics and Astronautics, AIAA Student Member.}, L{\'e}opold Cambier\footnote{Ph.D. Candidate, Institute for Computational and Mathematical Engineering.}, Juan J.~Alonso\footnote{Vance D. and Arlene C. Coffman Professor, Department of Aeronautics and Astronautics, AIAA Associate Fellow.} and Eric Darve\footnote{Professor, Department of Mechanical Engineering.}}
\affil{Stanford University, Stanford, CA, 94305}

\begin{document}

\maketitle

\begin{abstract}
Development of highly scalable and robust algorithms for large-scale CFD simulations has been identified as one of the key ingredients to achieve NASA's CFD Vision 2030 goals. In order to improve simulation capability and to effectively leverage new high-performance computing hardware, the most computationally intensive parts of CFD solution algorithms ---namely, linear solvers and preconditioners--- need to achieve asymptotic behavior on massively parallel and heterogeneous architectures and preserve convergence rates as the meshes are refined further. In this work, we present a scalable high-order implicit Discontinuous Galerkin solver from the SU2 framework using a promising preconditioning technique based on algebraic sparsified nested dissection algorithm with low-rank approximations, and communication-avoiding Krylov subspace methods to enable scalability with very large processor counts. The overall approach is tested on a canonical 2D NACA0012 test case of increasing size to demonstrate its scalability on multiple processing cores. Both the preconditioner and the linear solver are shown to exhibit near-linear weak scaling up to 2,048 cores with no significant degradation of the convergence rate.
\end{abstract}

\section{Nomenclature}

{\renewcommand\arraystretch{1.0}
\noindent\begin{longtable*}{l @{\quad=\quad} p{6.3cm}l @{\quad=\quad} p{6.3cm}}
$A$ & Matrix of linear system & $R^n$ & Spatial residual evaluated at timestep $n$\\
$D_k$ & Discrete element $k$ & $\frac{\partial R}{\partial u}$ & Spatial Jacobian\\
$F$ &  Convective flux & $s$ & Step size in $s$-step GMRES \\
$F^*$ &  Numerical flux & $t$  & Temporal coordinates\\
$H$ & Upper Hessenberg matrix & $U$ & Upper triangular matrix\\
$I$ & Identity matrix & $u$  & Conserved variables\\
$K$ & Number of discrete elements & $\hat{u}$  & Modal coefficients\\
$L$ & Lower triangular matrix & $\bar{u}$  & Nodal coefficients\\
$\ell$ & level in Nested Dissection ordering & $\vec{x}$ & Cartesian coordinates\\
$l_i$ & Lagrange polynomial&  $\alpha$ & Angle of attack\\
$M$ & Mass matrix&  $\Delta \tau$  & Time step\\
$M_\infty$ & Free-stream Mach number&  $\varepsilon$ & tolerance in spaND\\
$N$ & Dimension of Jacobian matrix & $\lambda$ & Spectral radius of the convective operator\\
$N_p$ & Number of degrees of freedom& $\rho$ & Density\\
$n_i$  & Normal vector at interfaces& $\psi_i$ & Legendre polynomial\\
$p$ & Polynomial order& $\Omega$& Physical domain\\
$Q$ & Orthogonal matrix & $\partial\Omega$& Boundary of physical domain\\
\end{longtable*}}

\section{Introduction}
\label{sec:intro}

\lettrine{A}{dvancements} in computational capabilities have enabled many large-scale simulations of complex fluid flow problems. With the expected arrival of next-generation exascale supercomputers in 2021 \cite{exasupercomputer}, highly efficient, robust, and scalable CFD solvers are needed to fully harness the power of High Performance Computing (HPC). For various HPC codes solving the Navier-Stokes equations and their various approximations, the most computationally-intensive parts are often contained in numerical algorithms such as linear solvers and preconditioners. As processor count and problem size increase, the time spent in these numerical algorithms, relative to the rest of the application, grows and quickly dominates the total execution time \cite{ang2012report}. This is especially true for high-order methods such as Discontinuous Galerkin (DG) methods that are computationally expensive. Hence, developing linear solvers and preconditioners suited to massive CFD problems and very large processor counts plays a crucial part in building scalable CFD solvers.

Unfortunately, despite decades of work in the area, the development of preconditioners for linear systems are still considered an art and the details depend on the specific problem being solved. Many options are available but none work well in all cases. Current options for large-scale implicit discretizations of the Navier-Stokes equations are limited to multigrid approaches, which are difficult to implement and yield limited benefits \cite{baker2012scaling}, and incomplete LU factorizations (ILU) which work for small problems but deteriorate rapidly for larger grids \cite{saad1994ilut}. In addition, many numerical algorithms are found to incur overwhelming communication costs when moving data between levels of the memory hierarchy or between processors over a network on massively-parallel computers \cite{hoemmen2010communication}, especially when synchronization points across the entire parallel computer are an integral part of the algorithm. Achieving scalable linear solvers thus requires not only robust and reliable preconditioning techniques, but also a dramatic shift in algorithm design with a focus on reducing communication.

In this paper, we present our recent efforts towards constructing a scalable flow solver using an implicit high-order DG discretization. The work described in this paper focuses on the development of ideas for the core preconditioned linear solver step and does not go into much detail in the description of the DG discretization in the SU2 solver. For the linear system, we apply a promising algebraic preconditioner that is based on hierarchical matrices and low-rank approximations, and has shown scalability for many partial differential equation (PDE) problems~\cite{cambier2020algebraic}. To overcome communication overhead, we implement variants of a Krylov subspace method that avoids/hides communication to solve the preconditioned linear system. A number of scalability studies are performed to demonstrate the scalability of the overall approach. The paper is organized as follows: Section \ref{sec:methodology} gives technical details of the different aspects of the proposed scalable CFD solver, Section \ref{sec:results} outlines the numerical results for the integrated approach, and Section \ref{sec:conclusion} concludes and describes the intended future work.

\section{Methodology}
\label{sec:methodology}
In this section, we present various components of the overall scalable CFD solver framework: a high-order DG flow solver, a novel preconditioning technique, and a communication-avoiding iterative linear solver. In combination, these methods can achieve truly scalable solution methods for complex CFD problems.
\subsection{Flow Solver}\label{sec:flowsolver}
\subsubsection{Governing equations}

In this work, we are interested in compressible form of the inviscid Euler equations of gas dynamics,
\begin{equation}
    \frac{\partial u}{\partial t} + \nabla \cdot F = 0.
\end{equation}

The governing equations are solved in a physical domain $\Omega$ with boundary $\partial\Omega$.

\subsubsection{Implicit Discontinuous Galerkin discretization}


\indent DG methods form a class of numerical methods for solving PDE problems with high order of accuracy \cite{hesthaven2007nodal}. In a Discontinuous Galerkin Finite Element Discretization (DG-FEM), we divide the computational domain $\Omega$ into $K$ discrete elements $D_k$,
\begin{equation}
    \Omega = \bigcup_{k=1}^K D_k.
\end{equation}

The numerical solution within each element $D_k$ is constructed using nodal or modal representations as

\begin{equation}
    u(\vec{x}, t) = \sum_{i=1}^{N_p}\hat{u}_i(t)\psi_i(\vec{x}) = \sum_{j=1}^{N_p}\bar{u}(x_j,t)l_j(\vec{x}).
\end{equation}

The number of degrees of freedom $N_p$ is defined by the order of the polynomial basis $p$. By allowing discontinuity at element boundaries, the weak form of the governing equations over each spatial element $D_i$ can be expressed using a nodal representation as

\begin{equation}
    \int_{D_k}\frac{\partial{\bar{u}}}{\partial{t}}l_jdV-\int_{D_k}F_i\frac{\partial{l_j}}{\partial{x_i}}dV+\oint_{\partial{D_k}}F_i^*n_il_jdA=0,\qquad j=1,...,N_p.
\end{equation}

A common choice of the numerical flux $F^*$ is the Roe flux \cite{roe1981approximate}. The implicit discretization with a backwards Euler temporal discretization leads to a linear system at every timestep $n$

\begin{equation}
    \left(\frac{1}{\Delta\tau^n}M+\frac{\partial{R}}{\partial{u}}\bigg\rvert^n\right)\Delta{u}^{n+1}=-R^n,\label{eq:linearsystem}
\end{equation}

where

\begin{equation}
    M = \int_{D_k}l_il_jdV,
\end{equation}
\begin{equation}
    R^n = -\int_{D_k}F_i(u^n)\frac{\partial{l_j}}{\partial{x_i}}dV+\oint_{\partial{D_k}}F^*_i(u^n)n_il_jdA.
\end{equation}

For stability, the timestep $\Delta\tau$ is chosen based on the CFL number

\begin{equation}
    \Delta\tau = \min(\Delta\tau_k) = \min\left(\frac{\mbox{CFL}}{\lambda_k}\right).
\end{equation}

For steady-state problems, pseudo time-stepping is used to accelerate convergence where the CFL number at the $n$th timestep is computed by

\begin{equation}
    \mbox{CFL}^n = \min\left(\mbox{CFL}^0\frac{\left\lVert{R^0}\right\rVert}{\left\lVert{R^{n-1}}\right\rVert},\mbox{CFL}^\infty\right).
\end{equation}

DG methods are local methods in the sense that, by allowing discontinuity at element boundaries, the high-order solution within each discrete element depends only on its immediate neighbors. This has significant advantages in a distributed-memory architecture as each partition of the computational domain relies only on its interior elements and a single layer of halo elements at the boundaries of the partition to construct high-order solutions. Such partitions typically have a low surface-to-volume ratio that makes the solver highly parallelizable. In an implicit  discretization, the locality of DG methods often gives rise to Jacobians with a sparse, block matrix structure. The low surface-to-volume ratio translates to minimal point-to-point communication costs in matrix kernels that scale well on increasing number of processors.

\subsubsection{Implementation}

The SU2 software suite \cite{palacios2013stanford,palacios2014stanford,economon2015} is an open-source collection of software tools written in C++ and Python for performing multi-physics simulation and design. It is built specifically for the analysis of PDEs and PDE-constrained optimization problems on general unstructured meshes with state-of-the-art numerical methods. The DG-FEM solver is one of the solvers in SU2. It supports all standard elements in two and three dimensions up to arbitrary polynomial orders with the option of local $p$-refinement. Other features of the DG-FEM solver include treatment of curved elements, treatment of viscous terms, shock-capturing capability, ADER-DG discretization, etc \cite{choi2019simple}. For the implicit DG discretization in SU2, the spatial Jacobian $\frac{\partial{R}}{\partial{u}}$ in Eq. \eqref{eq:linearsystem} is evaluated exactly using the automatic differentiation (AD) tool CodiPack \cite{sagebaum2017high} which allows the implicit formulation to leverage all features available in SU2.
\subsection{Sparsified Nested Dissection and TaskTorrent}

\subsubsection{Sparsified Nested Dissection algorithm}

Sparsified Nested Dissection (spaND) is a fast multilevel algorithm
for solving large sparse linear systems \cite{cambier2020algebraic}. Let $Ax=b$ be the linear system from Eq. \eqref{eq:linearsystem}. The algorithm first computes a Nested Dissection (ND) ordering of the matrix $A$. This defines interiors, separators, and interfaces at each dissection level $0 \leq \ell < \ell_{\max}$. Level $0$ is the leaf level and $\ell_{\max}$ is the top level.
At each level, interiors are separated by the ND separators and eliminating an interior does not create fill-in beyond its adjacent separator. Interfaces are defined as subsets of separators adjacent to a given pair of interiors on each side of the separator.

The algorithm then proceeds level by level, from the leaf level to the top level.
At each level $\ell$:
\begin{itemize}
    \item Interiors at level $\ell$ are eliminated using row-pivoted block LU factorization.  This is the same algorithm as any sparse direct method. With $A_{ss} = PLU$ we find
    \begin{equation}
    \begin{bmatrix} L^{-1} P^\top & \\ -A_{ns} A_{ss}^{-1} & I \end{bmatrix} \begin{bmatrix} A_{ss} & A_{sn} \\ A_{ns} & A_{nn} \end{bmatrix} 
    \begin{bmatrix} U^{-1} & -A_{ss}^{-1} A_{sn} \\ & I \end{bmatrix} = \begin{bmatrix} I & \\ & A_{nn} - A_{ns} A_{ss}^{-1} A_{sn} \end{bmatrix}.
    \end{equation}
    This step introduces fill-in on $A_{nn}$ because of the $A_{ns} A_{ss}^{-1} A_{sn}$ term. Notice that many rows (respective columns) in $A_{ns}$ (resp. $A_{sn}$) are zero and, as such, do not have to be updated.
    \item Interfaces at level $\ell$ are scaled using a row-pivoted block LU factorization. If $A_{ss} = PLU$, we have
    \begin{equation}
    \begin{bmatrix} L^{-1} P^\top & \\ & I \end{bmatrix} 
    \begin{bmatrix} A_{ss} & A_{sn} \\ A_{ns} & A_{nn} \end{bmatrix} 
    \begin{bmatrix} U^{-1} & \\ & I \end{bmatrix} = 
    \begin{bmatrix} I & L^{-1} P^\top A_{sn} \\ A_{ns} U^{-1} & A_{nn} \end{bmatrix}.
    \end{equation}
    Scaling is required for accuracy as scaling the diagonal blocks leads to a much lower number of iterative method steps \cite{cambier2020algebraic}. Note that this is not an approximation and does not create any fill-ins. It does, however, balance the matrix since all diagonal blocks now have the same unit norm. Note that in this step, it is preferable to balance $L$ and $U$ so that $\|U^{-1}\| \approx \|L^{-1}\|$. We do so by splitting the diagonal of the upper-triangular matrix from the LU factorization evenly between $L$ and $U$.
    \item Interfaces at level $\ell$ are sparsified using a low-rank approximation (in practice, rank-revealing QR). 
    Let $s$ be an interface and $n$ all the neighbors of $s$ in the trailing matrix $A$. 
    We first compute $Q_s = \begin{bmatrix} Q_{sc} & Q_{sf} \end{bmatrix}$ such that
    \begin{equation}
        \begin{bmatrix} A_{sn} & A_{ns}^\top \end{bmatrix} = Q_{sc} \begin{bmatrix} W_{cn} & W_{nc}^\top \end{bmatrix} + Q_{sf} \begin{bmatrix} W_{fn} & W_{nf}^\top \end{bmatrix},
    \end{equation}
    with $\left\| \begin{bmatrix} W_{fn} & W_{nf}^\top \end{bmatrix} \right\|_2 = \mathcal{O}(\varepsilon)$
    where $\varepsilon$ is a user-prescribed tolerance. 
    Let $w$ be the remaining degrees of freedom disconnected from $s$.
    Given this low-rank approximation, the trailing matrix can be factorized as
    \begin{equation}
        \begin{bmatrix} Q^\top  &   & \\
                                & I & \\
                                &   & I \end{bmatrix}
        \begin{bmatrix} I      & A_{sn} &        \\
                        A_{ns} & A_{nn} & A_{nw} \\
                               & A_{wn} & A_{ww} \end{bmatrix}
        \begin{bmatrix} Q   &   & \\
                            & I & \\
                            &   & I \end{bmatrix}
        =
        \begin{bmatrix} I       &               & W_{cn}        &           \\
                                & I             & \varepsilon   &           \\
                        W_{nc}  & \varepsilon   & A_{nn}        & A_{nw}    \\
                                &               & A_{wn}        & A_{ww}    \end{bmatrix}.
    \end{equation} 
    Assuming $\varepsilon \approx 0$, the rows and columns corresponding to the variable $f$ are effectively 0, so $f$ is approximately eliminated. In addition, this procedure did not introduce any fill-in on the neighbors $n$ of $s$ ($A_{nn}$ is unchanged). So this (approximately) eliminated some of the non-interior unknowns without introducing any fill-in in the trailing matrix.
    
    In practice, we first compute all $Q_i$ without updating the trailing matrix. Only then do we compress every $A_{ij}$ block and replace it with $A_{ij}^+ = Q_{ic}^\top A_{ij} Q_{jc}$. This leads to the same approximation as described above, but is more concurrent since every rank-revealing QR can be done simultaneously.
    \item We then merge all the clusters and proceed to the next level.
\end{itemize}  

At the end, the algorithm produces an approximate factorization of $A$, $A \approx \prod_{i} F_i$ where $F_i$ is a sparse triangular matrix (from the elimination or block scaling) or a sparse orthogonal matrix (from the sparsification).
This is then used as a preconditioner for Krylov iterative methods, such as Conjugate Gradient \cite{hestenes1952methods}, Generalized Minimum RESidual (GMRES) \cite{saad1986gmres}, etc.

\begin{figure}[H]
    \centering
    \subfloat[][\label{sfig:leaf}Initial matrix. Each dot corresponds to four unknowns in the system $Ax=b$]{
        \includegraphics[trim={0 0 0 1cm},clip,width=0.4\textwidth]{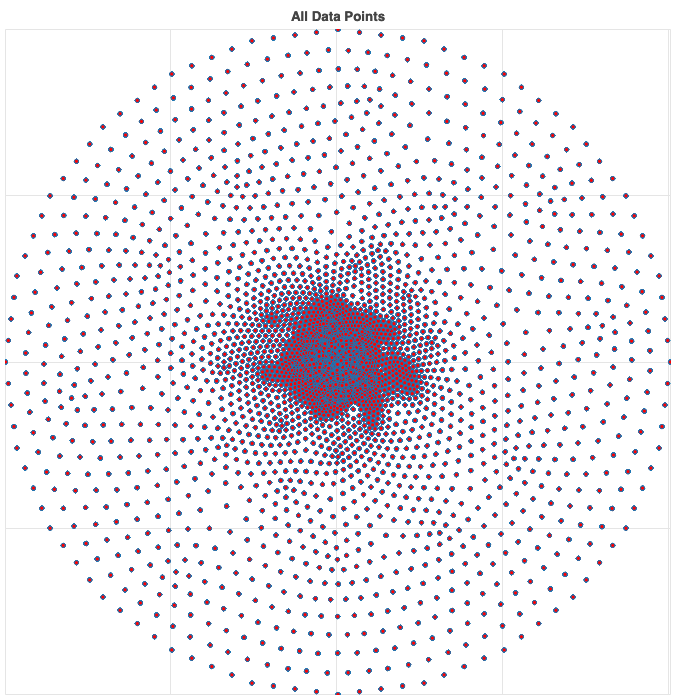}
    } \;\;
    \subfloat[][\label{sfig:elim}After interiors ($\ell=5$) elimination]{
        \includegraphics[trim={0 0 0 1cm},clip,width=0.4\textwidth]{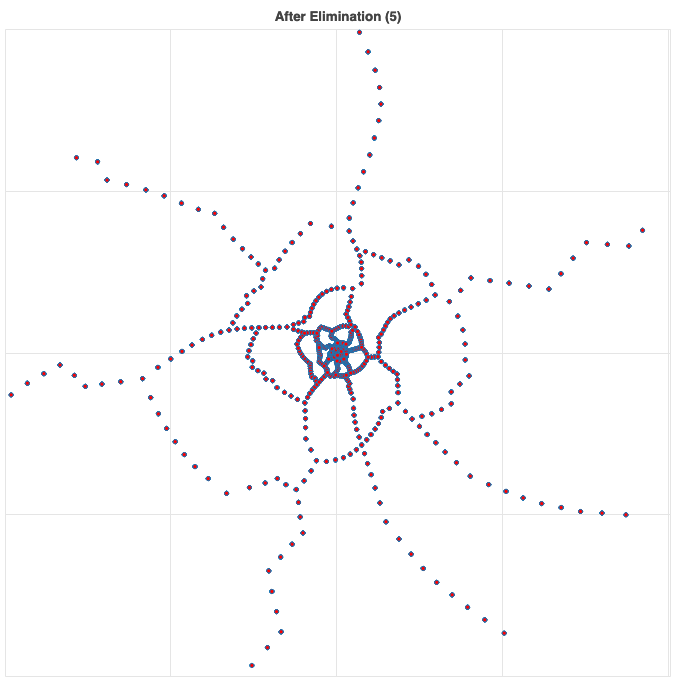}
    } \\
    \subfloat[][\label{sfig:spars}After interface ($\ell=5$) sparsification. All interfaces are reduced in size.]{
        \includegraphics[trim={0 0 0 1cm},clip,width=0.4\textwidth]{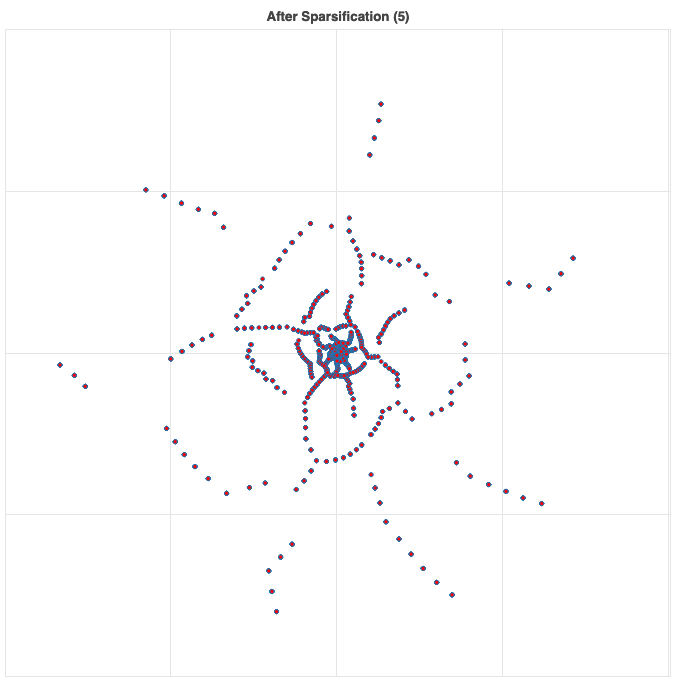}
    } \;\;
    \subfloat[][\label{sfig:top}Last top separator. A typical ND separator would be a line cutting through the entire domain.]{
        \includegraphics[trim={0 0 0 1cm},clip,width=0.4\textwidth]{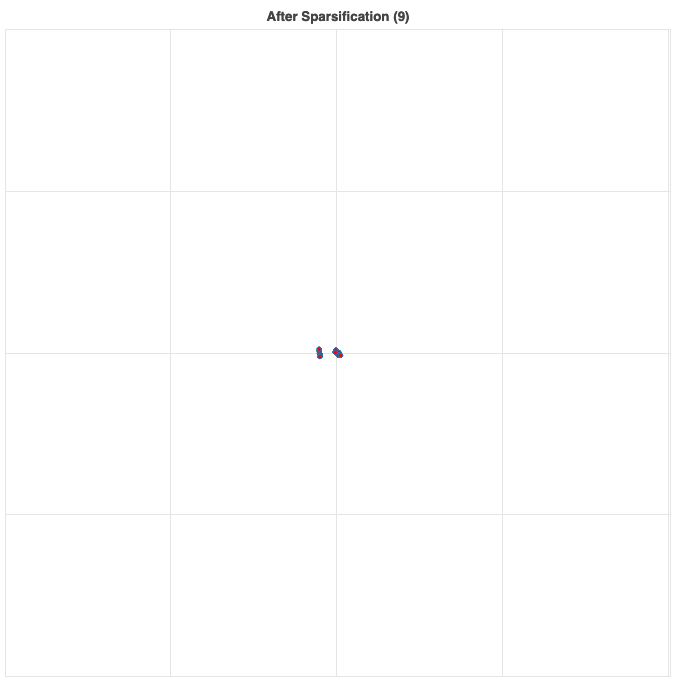}
    }
    \caption{Illustration of the spaND algorithm. At every level, interiors are eliminated and the remaining interfaces are sparsified. Interfaces are then merged and the algorithm proceeds to the next level. This leads to separators of size $\mathcal{O}(1)$ instead of lines of size $\mathcal{O}(N^{1/2})$.}
    \label{fig:elim_spars}
\end{figure}

\Cref{fig:elim_spars} illustrates the algorithm applied to a linear system generated from implicit discretization of a typical 2D airfoil mesh.
\Cref{sfig:leaf} shows all degrees of freedom in the system.
\Cref{sfig:elim} shows the remaining degrees of freedom after 5 levels of ND elimination. We clearly see the eliminated interiors (the large white areas) and the remaining separators and interfaces. \Cref{sfig:spars} shows the effect of interface sparsification. All interfaces (i.e., subsets of a separator separating a pair of interiors) are sparsified, which reduces their size without introducing fill-in. This reduces the sizes of all separators in the systems.
\Cref{sfig:top} shows the top separator before its final elimination. 
Instead of a usual ND separator cutting across the domain, only a few points are left.

By repeating this process at every level, this algorithm keeps the separator size small.
In 2D problems, separators typically have size $\mathcal{O}(1)$ (instead of $\mathcal{O}(N^{1/2})$ without sparsification), while in 3D, separators now have size $\mathcal{O}(N^{1/3})$ (instead of $\mathcal{O}(N^{2/3})$). 
This leads to a great reduction in computational cost.
A direct method using an ND ordering usually has complexity $\mathcal{O}(N^{3/2})$ for 2D problems. 
In contrast, assuming the separators shrink to size $\mathcal{O}(1)$, spaND has complexity $\mathcal{O}(N)$ in 2D, for both factorization and solve time. It is then typically used as a preconditioner coupled with CG or GMRES.
Assuming a small and slowly growing iteration count, this leads to a linear or near-linear time algorithm.
\subsubsection{Parallel implementation using TaskTorrent}

We finally consider a parallel version of the spaND algorithm. 
To do so, we use a task-based runtime system in C++, TaskTorrent \cite{cambier2020tasktorrent}.
TaskTorrent is a lightweight and distributed task-based runtime system where computations are expressed as a directed acyclic graph (DAG) of tasks using a parametrized task graph formulation. Tasks in the DAG are then run as soon as ready by the runtime. Active messages (a form of one-sided asynchronous communications combining a function and a payload) let tasks trigger other tasks on remote nodes. TaskTorrent uses C++ threads for intra-node parallelism and MPI for inter-nodes communications. 

The key step in parallelizing spaND using TaskTorrent is the definitions of the DAG of tasks.
We do so by transforming every block operation into a corresponding task.
This is done at the granularity of interfaces. 
Furthermore, we parallelize each level independently. This means there is one distinct DAG of task per level in the algorithm.
In TaskTorrent, this requires identifying the number of incoming dependencies for each task, the computational routine of the task itself, and the outgoing dependencies of each task.

Consider for instance the rank-revealing QR factorization related to the sparsification of interface $s$, denoted \texttt{geqp3(s)}. 
Let $n_1, \dots, n_k$ denote all the neighboring interfaces of $s$.
\begin{itemize}
    \item The task requires $A_{sn_1}, \dots, A_{sn_k}, A_{n_1s}, \dots, A_{n_ks}$. As such, \texttt{geqp3(s)} has $2k$ incoming dependencies.
    \item Given $A_{sn_1}, \dots, A_{sn_k}, A_{n_1s}, \dots, A_{n_ks}$, the task computes $\begin{bmatrix} Q_c & Q_f \end{bmatrix} = \text{geqp3}\left( \begin{bmatrix} A_{sn} & A_{ns}^\top \end{bmatrix}\right)$.
    This is done using LAPACK's \texttt{geqp3} function.
    \item $Q_c$ is used during the compression steps on blocks $A_{sn_1}, \dots, A_{sn_k}, A_{n_1s}, \dots, A_{n_ks}$. The compression of block $A_{ij}$ is denoted \texttt{ormqr(i,j)} and performs $A_{ij} \leftarrow Q_{ic}^\top A_{ij} Q_{jc}$.
    As such, \texttt{geqp3(s)} has to fulfill the dependencies of tasks \texttt{ormqr(s,n1)}, ..., \texttt{ormqr(s,nk)}, \texttt{ormqr(n1,s)}, ..., \texttt{ormqr(nk,s)} and send the basis $Q_c$ to the corresponding MPI ranks.
\end{itemize}

\begin{figure}
    \centering
    \includegraphics[width=0.3\textwidth]{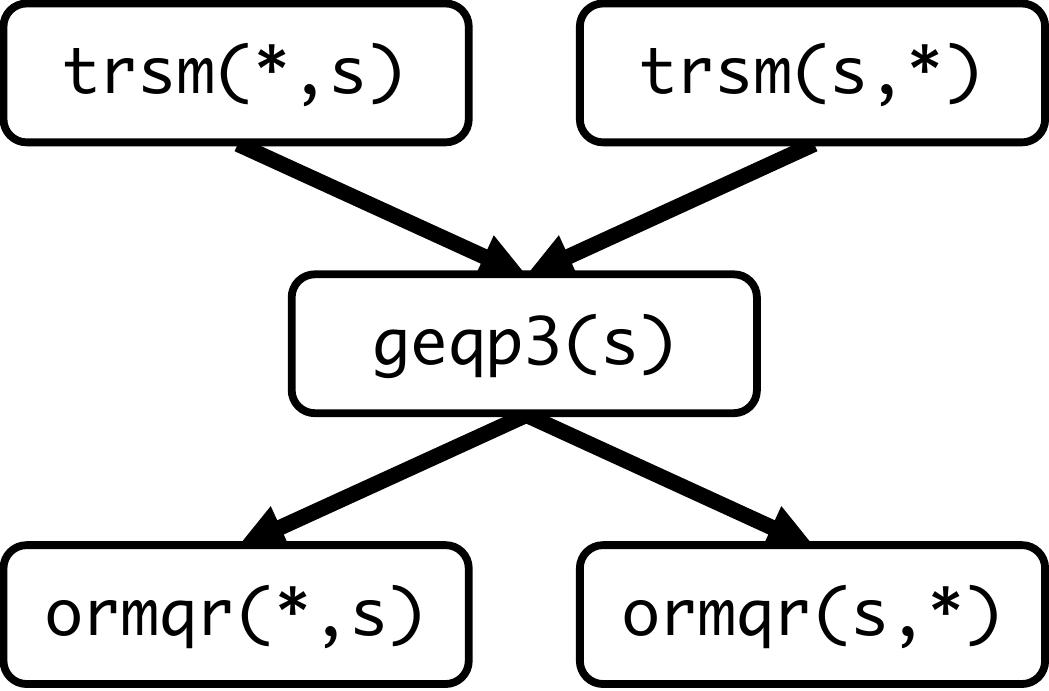}
    \includegraphics[width=0.2\textwidth]{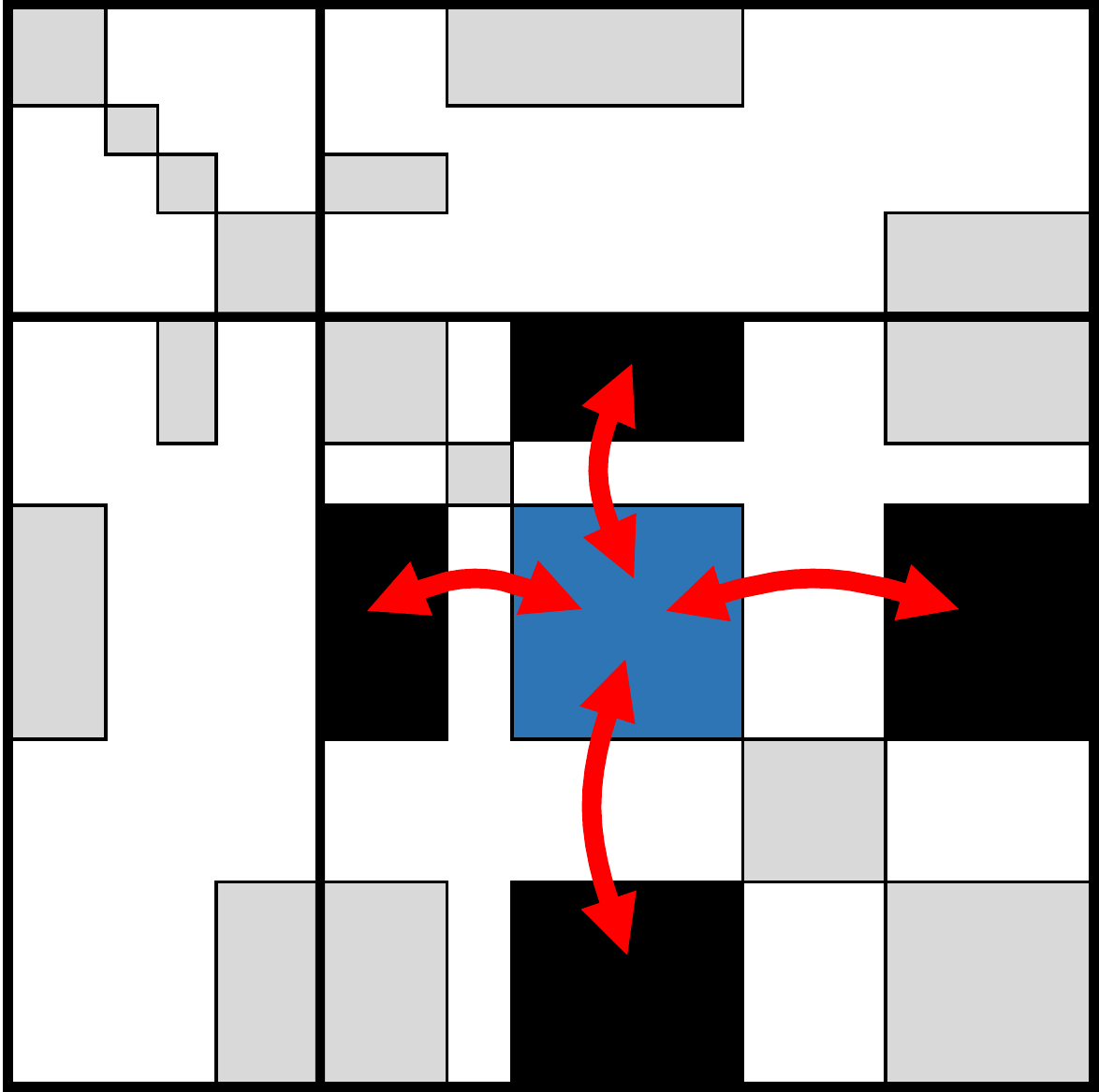}
    \caption{Task-based sparsification of an interface $s$. The left shows the local DAG, with the in- and out-dependencies of task \texttt{geqp3(s)}. The right figure shows the trailing matrix, with the block $(s,s)$ in blue. The first four diagonal blocks are interiors (at this point, all eliminated) and other diagonal blocks are interfaces. Red arrows show incoming and outgoing dependencies of task \texttt{geqp3(s)}.}
    \label{fig-pspand:sparsification}
\end{figure}

Fig. \ref{fig-pspand:sparsification} illustrates the related DAG and the dependencies in the trailing matrix. All tasks in spaND are parallelized using TaskTorrent in a similar way to miminize load imbalance. The overall framework of spaND and TaskTorrent are incorporated in the SU2 codebase.

\subsection{Communication-avoiding (CA) Krylov Subspace methods}\label{seckrylov}
Krylov subspace methods are a class of iterative methods commonly applied to large, sparse linear systems in CFD codes. However, on large machines, traditional Krylov solvers do not scale well due to their reliance on expensive communication steps in various kernels, particularly those involving global communication constructs. To alleviate performance bottlenecks due to communication overhead, we implement CA Krylov subspace methods to achieve asymptotic performance improvements when solving large-scale linear systems.

In this work, we are interested in the GMRES solver and its CA variants. The classical GMRES algorithm shown in Algorithm \ref{alg1} will incur Point-to-Point communication costs during Sparse Matrix-Vector (SpMV) multiplication (line 3) in a distributed memory setting. In addition, dot product computations (lines 5, 8) in the Modified Gram-Schmidt (MGS) algorithm require global reductions which are very costly synchronization steps, especially in the presence of system noise and load imbalances. The number of global reductions in MGS also grows quadratically with the iteration count. These kernels make the performance of the GMRES solver communication-bound.

To improve the performance of GMRES, CA variants of GMRES have been proposed \cite{hoemmen2010communication,ghysels2013hiding,yamazaki2017improving}. This work employs the $s$-step CAGMRES algorithm in \cite{hoemmen2010communication}, which attempts to avoid communication by generating $s$ Krylov basis vectors at once and orthogonalize them using a block version of the Classical Gram-Schmidt (CGS) algorithm followed by a Cholesky QR factorization (CholQR). This reduces the number of reductions required to one per $s$ vectors, achieving a reduction in the communication latency by a factor of $s$ at the expense of additional arithmetic operations. However, at large values of $s$, consecutive SpMVs produce vectors that converge to the principal eigenvector of the matrix, which then introduce numerical instability in subsequent orthogonalization steps. To overcome this, a different basis function is often incorporated to condition the $s$ block. Re-orthogonalization schemes, such as two CGSs and two CholQRs, that do not require extra communication are also introduced to minimize the orthogonality error \cite{yamazaki2020low, swirydowicz2020low}. A skeleton of the $s$-step CAGMRES algorithm with single-reduce re-orthogonalization is shown in Algorithm \ref{alg2}. Only one single reduction is required for each $s$ vector (line 6). This is achieved by lagging the re-orthogonalization and re-normalization of previous $s$ vectors (line 7, 8) and combining them with the orthogonalization and normalization of $s$ vectors in the next iteration (line 12, 13). A detailed algorithm can be found in \cite{yamazaki2020low}. It is worth pointing out that $s$ times SpMV in line 4 can also be applied in a CA manner, see \cite{hoemmen2010communication}. However, it is challenging to apply preconditioners in a similar CA fashion. Therefore, the CA variant of SpMV kernel is left for future work. Similar to spaND, CAGMRES algorithm in this work is implemented in the SU2 codebase.

\begin{algorithm}[ht!]
\caption{Classical GMRES} \label{alg1}
\hspace*{\algorithmicindent} \textbf{Input:} $n\times n$ matrix $A$, right hand side vector $b$, initial guess vector $x_0$ \\
\hspace*{\algorithmicindent} \textbf{Output:}{ $x$, solution to the linear system $Ax=b$.} 
\begin{algorithmic}[1]
\State $r := b - Ax_0, q_0 := r/||r||_2$
\For{$j = 0,1,...$}
\State $v = Aq_j$ \Comment{Sparse Matrix-Vector Multiplication (SpMV)}
\For{$i = 0,...,j$}\Comment{Modified Gram-Schmidt (MGS)}
\State $h_{ij} = v^Tq_i$ \Comment{global reduction}
\State  $v = v - h_{ij}q_i$
\EndFor\label{euclidendwhile}
\State$h_{j+1,j} = ||v||_2$ \Comment{global reduction}
\State$q_{j+1} = v/h_{j+1,j}$
\State Apply Givens rotation to update matrix $H_j$
\State Check for convergence
\EndFor
\State$y = $ argmin$||(H_jy-||r||_2e_0)||_2$
\State$x = x_0 + Q_j y$
\end{algorithmic}
\end{algorithm}

\begin{algorithm}[ht!]
\caption{$s$-step CAGMRES with single-reduce re-orthogonalization}  \label{alg2}
\hspace*{\algorithmicindent} \textbf{Input:} $n\times n$ matrix $A$, right hand side vector $b$, initial guess vector $x_0$, step size \textit{s}. \\
\hspace*{\algorithmicindent} \textbf{Output:}{ $x$, solution to the linear system $Ax=b$.} 
\begin{algorithmic}[1]
\State $r := b - Ax_0, q_0 := r/||r||_2$
\For{$j = 0,s...$}
\For{$k = j...j+s-1$}
\State $Q_{:,k+1} = AQ_{:,k}$ \Comment{\textit{s} times SpMV}
\EndFor

\State $[R_{:,j-s:j-1}, R_{:,j:j+s}] = Q_{:, 0:j+s}^T[Q_{:, j-s:j-1}, Q_{:,j:j+s}]$ \Comment{single reduction}
\State CGS($Q_{:,0:j-s-1}, Q_{:,j-s:j-1}, R_{:,j-s:j-1}, R_{:,j:j+s}$) \Comment{Re-orthogonalize $Q_{:,j-s:j-1}$}
\State CholQR($Q_{:,j-s:j-1}, R_{:,j-s:j-1}, R_{:,j:j+s}$) \Comment{Re-normalize $Q_{:,j-s:j-1}$}
\State Assemble upper Hessenberg matrix $H_{j-1}$
\State Apply Givens rotation to update matrix $H_{j-1}$
\State Check for convergence
\State CGS($Q_{:,0:j-1}, Q_{:,j:j+s}, R_{:,j:j+s}$)  \Comment{Orthogonalize $Q_{:,j:j+s}$}
\State CholQR($Q_{:,j:j+s}, R_{:,j:j+s}$) \Comment{Normalize $Q_{:,j:j+s}$}
\EndFor
\State$y = $ argmin$||(H_{j-1}y-||r||_2e_0)||_2$
\State$x = x_0 + Q_{:,0:j-1} y$
\end{algorithmic}
\end{algorithm}

\section{Numerical Results}
\label{sec:results}

\subsection{Simulation Setup}

To demonstrate our scalable hierarchical CFD solver, we perform numerical experiments on a 2D NACA0012 airfoil in an inviscid compressible flow at an angle of attack $\alpha = 2^{\circ}$ and a free-stream Mach number of $M_\infty = 0.3$. 
The governing equations are discretized using the implicit DG formulation described in Section \ref{sec:flowsolver}. The computational mesh is a structured O-mesh with $p=1$ quadrilateral elements. The flow field is initialized uniformly with free-stream values and the simulations are converged until the density residual norm decreases by at least 8 orders of magnitude. As spaND is capable of factorizing ill-conditioned systems \cite{cambier2020algebraic}, the initial CFL number is boosted to $\mbox{CFL}^0 = 1000$ to accelerate convergence. 

\subsection{Scalability study of spaND}

\begin{figure}[h]
    \centering
    \begin{tikzpicture}
    \begin{axis}[
        legend cell align=left,
        axis lines = left,
        width=3in,
        height=2.5in,
        xlabel = Nonlinear iteration,
        ylabel = Relative residual ($\rho$),
        legend pos=outer north east,
        ytick={-10,-8,-6,-4,-2,0},
        yticklabels={$10^{-10}$,$10^{-8}$,$10^{-6}$,$10^{-4}$,$10^{-2}$,$10^0$},
    ]
    \addplot [black,mark=otimes]
    table{data/convergence_1.dat};
    \addplot [red,mark=*]
    table{data/convergence_2.dat};
    \addplot [blue,mark=diamond]
    table{data/convergence_3.dat};
    \addplot [hcgreen,mark=square]
    table{data/convergence_4.dat};
    \addplot [hcorange,mark=triangle]
    table{data/convergence_5.dat};
    \legend{Test 1,Test 2,Test 3,Test 4,Test 5};
    \end{axis}
    \end{tikzpicture}
    \caption{Relative convergence versus number of nonlinear iterations of Tests 1-5.}
    \label{fig:convergence}
\end{figure}
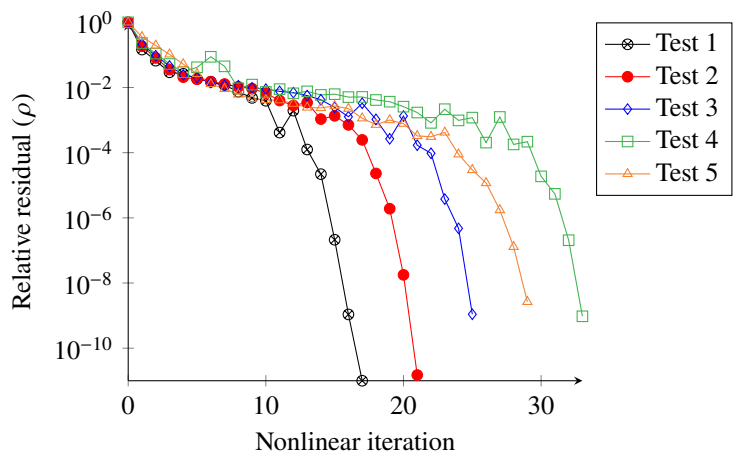

For scalability studies, we perform weak scaling analysis by varying the total number of elements in the mesh from 16k to 4M. As each first-order quadrilateral element has 4 nodal degrees of freedom and each node is represented by 4 conserved variables in 2D, the dimension of the resulting Jacobian matrix $N$ varies from 260k to 67M. The number of processing cores\footnote{Tests performed on a cluster equipped with dual-sockets and 16 cores Intel(R) Xeon(R) CPU E5-2670 0 @2.60GHz with 32GB of RAM per node} are scaled proportionally with respect to $N$. In this set of studies, a standard GMRES linear solver is used to solve the preconditioned linear system with a tolerance of $10^{-3}$ and maximum iteration count of $20$. The statistics of all tests as well as tunable parameters used in the overall solver setup are tabulated in Table \ref{table:1}. The convergence histories are shown in Fig. \ref{fig:convergence} and all tests converge within $35$ nonlinear iterations.

\begin{table}[!ht]
\caption{Summary of statistics and tunable parameters of tests}
\centering
\begin{tabular}{M{1cm}M{1.5cm}M{2cm}M{1.5cm}M{1.8cm}M{1cm}M{1.5cm}M{1.5cm}}
 \hline
 Test & No. of elements & Jacobian dimension, $N$ & No. of cores & SpaND tolerance, $\varepsilon$  & $\ell_{\max}$ & GMRES tolerance & GMRES max. steps \\
 \hline\hline
 1 & 16,384 & 262,144 & 8 & $10^{-3}$ & 10 & $10^{-3}$ & 20\\
 2 & 65,536 & 1,048,576 & 32 & $10^{-3}$ & 12 & $10^{-3}$ & 20\\
 3 & 262,144 & 4,194,304 & 128 & $10^{-3}$ & 14 & $10^{-3}$ & 20\\
 4 & 1,048,576 & 16,777,216 & 512 & $10^{-3}$ & 16 & $10^{-3}$ & 20\\
 5 & 4,194,304 & 67,108,864 & 2,048 & $10^{-3}$ & 18 & $10^{-3}$ & 20\\
 \hline
\end{tabular}
\label{table:1}
\end{table}

At each nonlinear iteration, spaND takes in the Jacobian matrix and generates an approximate factorization that is used as a preconditioner for linear solvers. 
Due to the nature of the discretization scheme, the Jacobian matrix is unsymmetric but with a symmetric sparsity pattern and a natural block structure. The mesh has a regular, structured topology such that each interior element is connected to four adjacent neighboring elements. However, due to the need to refine near-wall flow, the mesh is highly non-uniform with a fine mesh around the airfoil and coarse mesh near far-field boundaries.
During the pre-processing stage, we map the non-uniform mesh to a regular grid by assigning to each element a tuple $(i,j)$ of integers with $i$ increasing monotonically with the radius and $j$ with the angle, respectively.
The matrix is then partitioned using those $(i,j)$ coordinates with a standard recursive bisection algorithm. 

We then compute the factorization of the linear systems using spaND with TaskTorrent at every nonlinear step. 
We use partial pivoted LU as the block scaling algorithm and a tolerance of $\varepsilon = 10^{-3}$ for low-rank approximation.
Since the distribution of matrix ranks is not known beforehand, spaND distributes an equal number of columns of the matrix to each MPI rank at the beginning of the algorithm, using a 1D mapping of columns to MPI ranks. This indicates that if matrix ranks generated from low-rank approximation are higher in some parts of the domain than in others, load balancing may be sub-optimal as different amounts of computation is required in different processing cores.

\begin{figure}[h!]
    \centering
    \subfloat[$\ell = 4$. Ranks from 40 (blue) to 296 (brown), average of 162]{
        \includegraphics[width=0.4\textwidth]{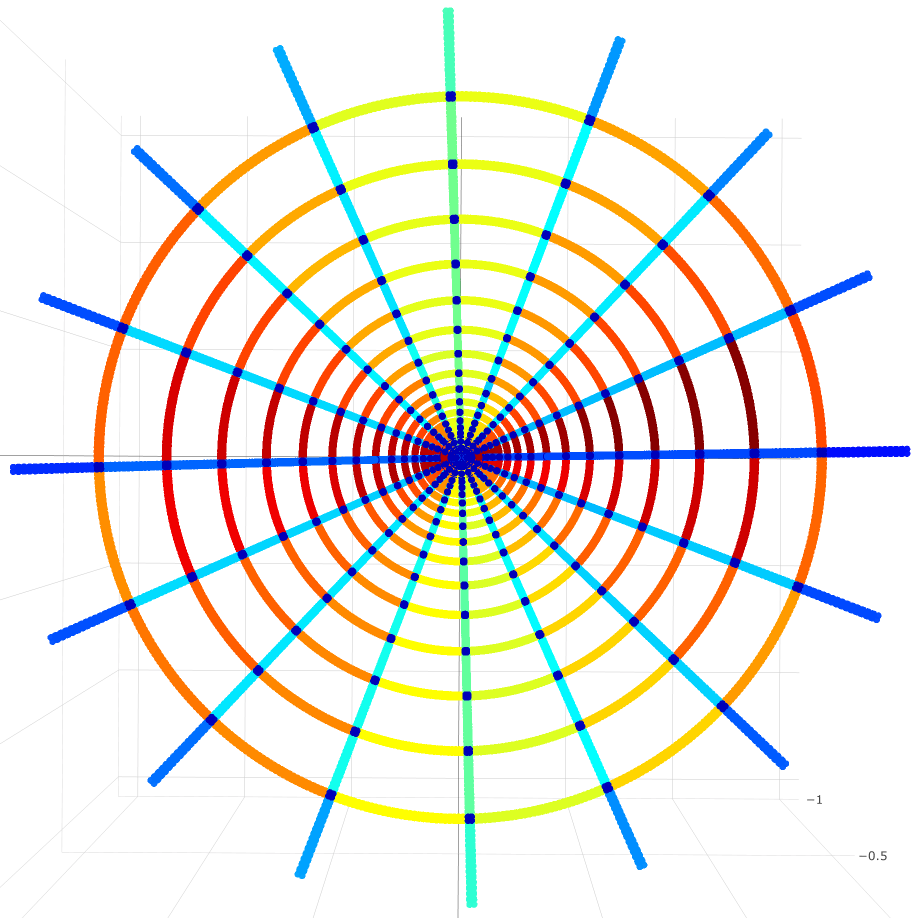}
    }\;
    \subfloat[$\ell = 6$. Ranks from 47 (blue) to 353 (brown), average of 198]{
        \includegraphics[width=0.4\textwidth]{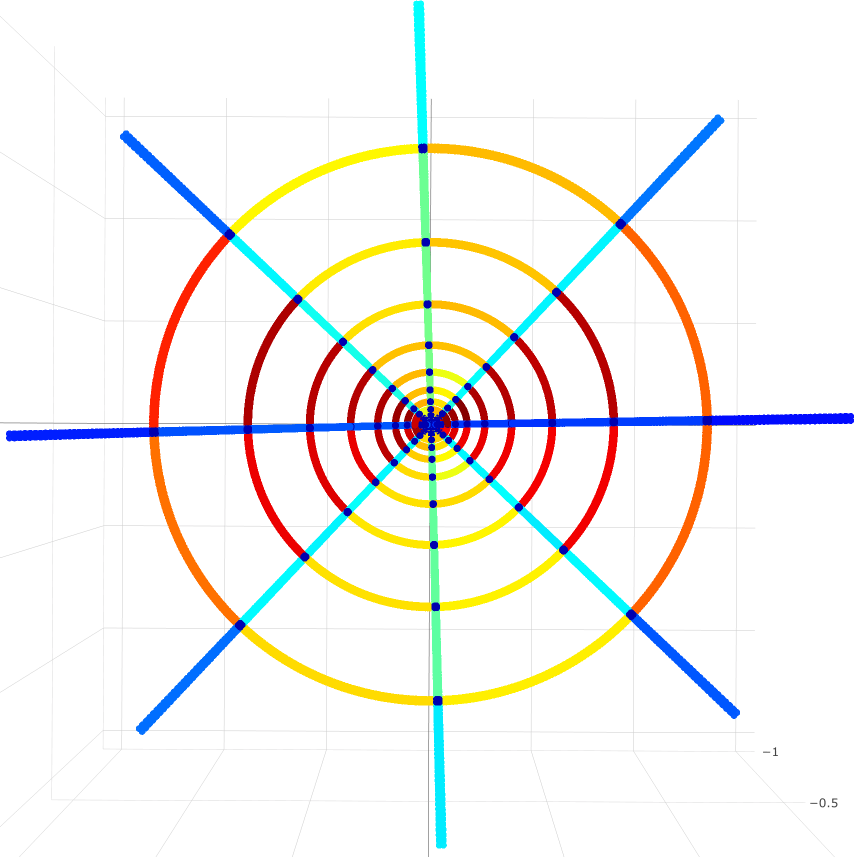}
    } \\
    \subfloat[$\ell = 6$, zoom on the airfoil.]{
        \includegraphics[width=0.4\textwidth]{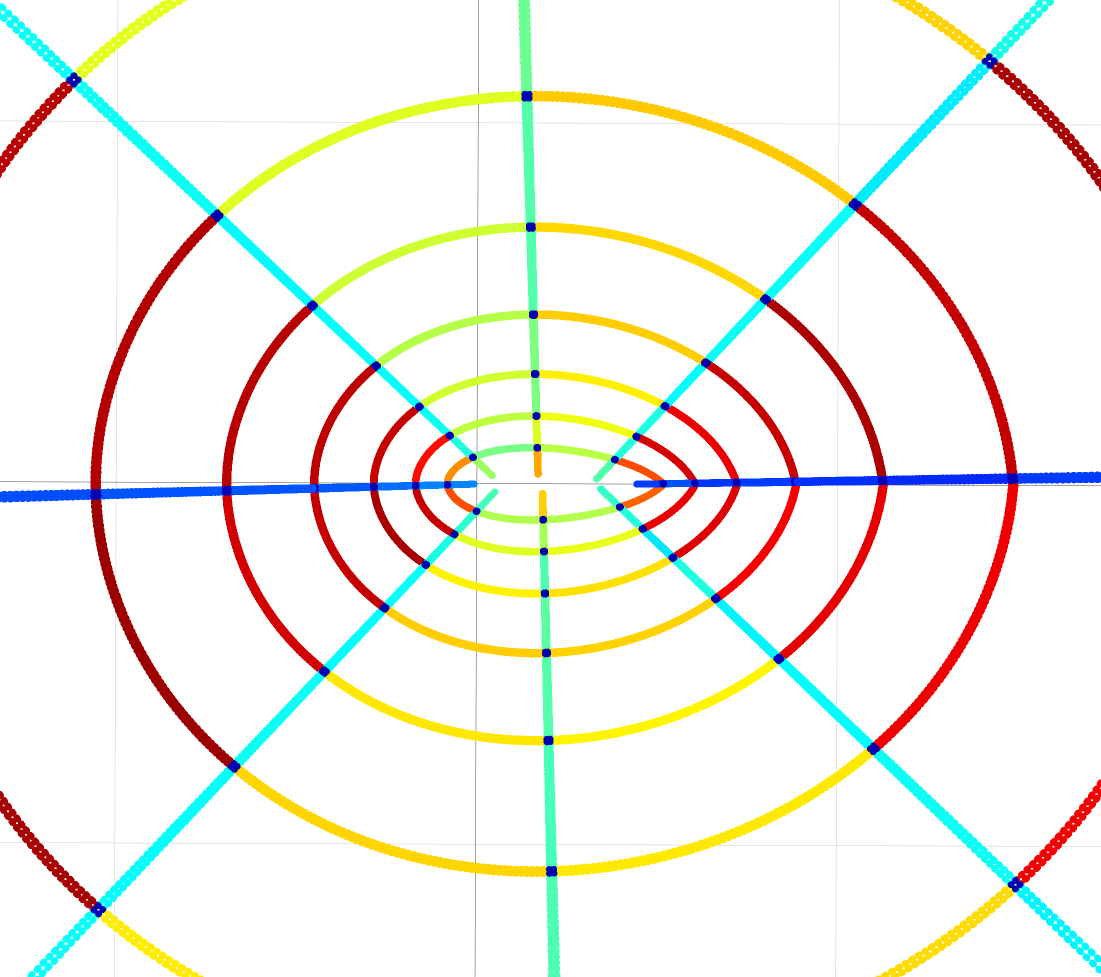}
    }\;
    \subfloat[$\ell = 8$. Ranks from 40 (blue) to 385 (brown), average of 231]{
        \includegraphics[width=0.4\textwidth]{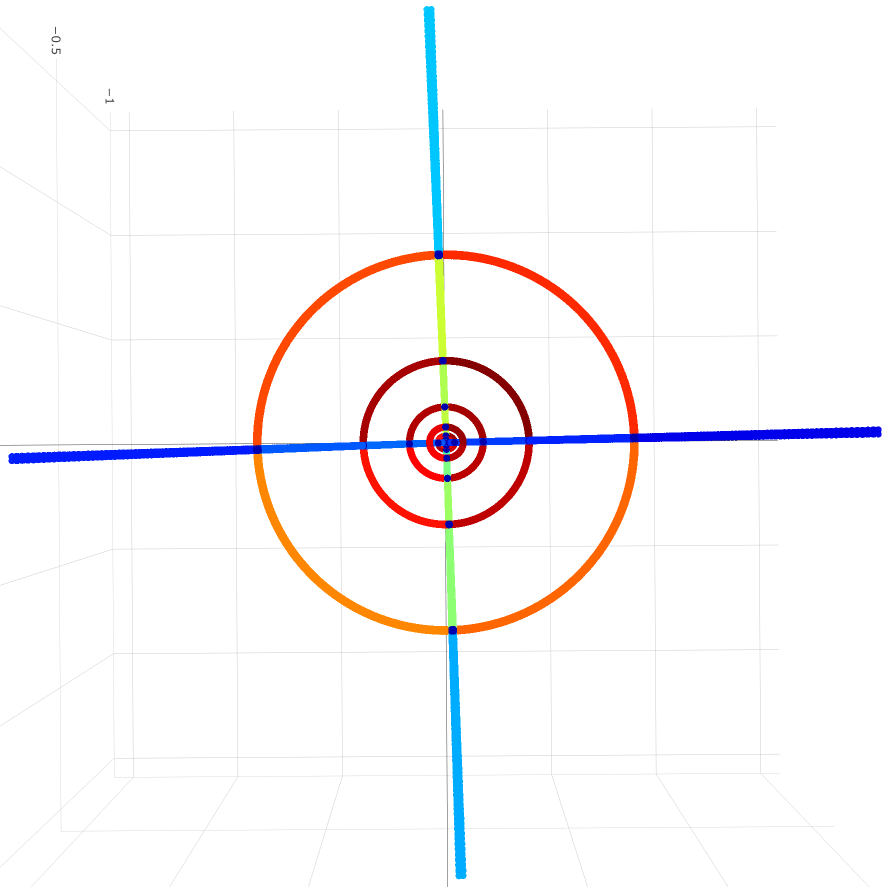}
    }
    \caption[Naca ranks]{Matrix rank distribution from rank-revealing QR at various level $\ell$ with $\ell_{\max}= 14$ in Test 3 at the last nonlinear iteration (iteration 25).}
    \label{fig-pspand:naca_ranks_geo_plots}
\end{figure}

To investigate this, we visualize the matrix rank distribution throughout the whole computational domain. Fig. \ref{fig-pspand:naca_ranks_geo_plots} shows a typical matrix rank distribution generated by spaND. We note that there is a strong directionality, with matrix ranks much higher in the flow direction, but smaller in the orthogonal direction. 
This presents some challenges since, as a consequence of wide range of matrix ranks, the load balancing is much less favorable.

\begin{figure}[h!]
    \centering
    \subfloat[Factorization time]{
        \begin{tikzpicture}\label{fact}
            \begin{axis}[
            axis lines = left,
                ymin=0,width=0.45\textwidth,xlabel={Nonlinear iteration},ylabel={Factorization time{, s}},
            ] 
                \addplot[black,mark=otimes]         table[x=step,y=tfact,discard if not={N}{262144}]   {data/naca_plu_all_tol1e-3_gmres1e-3.dat};
                \addplot[red,mark=*]                table[x=step,y=tfact,discard if not={N}{1048576}]  {data/naca_plu_all_tol1e-3_gmres1e-3.dat};
                \addplot[blue,mark=diamond]         table[x=step,y=tfact,discard if not={N}{4194304}]  {data/naca_plu_all_tol1e-3_gmres1e-3.dat};
                \addplot[hcgreen,mark=square]       table[x=step,y=tfact,discard if not={N}{16777216}] {data/naca_plu_all_tol1e-3_gmres1e-3.dat};
                \addplot[hcorange,mark=triangle]    table[x=step,y=tfact,discard if not={N}{67108864}] {data/naca_plu_all_tol1e-3_gmres1e-3.dat};
            \end{axis}
        \end{tikzpicture}
    }
    \subfloat[Ranks at the last nonlinear iteration. Average (solid line, diamond) and maximums (dashed line, circle).]{
        \begin{tikzpicture}\label{ranks}
            \begin{semilogyaxis}[
            axis lines = left,
                width=0.45\textwidth,
                ymin = 100,
                xlabel={Level from top},ylabel={Matrix Ranks},
                xtick={1, 6, 11, 16},xticklabels={15, 10, 5, 0},
                ytick={100,250,500,1000},yticklabels={100,250,500,1000},
                legend pos=outer north east,
                ]
                \addplot[mark=diamond,black]        table[x expr=\thisrow{lvl}+8,y=average] {data/naca_ranks_250k_tol1e-3_PLU.dat};
                \addplot[mark=diamond,red]          table[x expr=\thisrow{lvl}+6,y=average] {data/naca_ranks_1M_tol1e-3_PLU.dat};
                \addplot[mark=diamond,blue]         table[x expr=\thisrow{lvl}+4,y=average] {data/naca_ranks_4M_tol1e-3_PLU.dat};
                \addplot[mark=diamond,hcgreen]      table[x expr=\thisrow{lvl}+2,y=average] {data/naca_ranks_16M_tol1e-3_PLU.dat};
                \addplot[mark=diamond,hcorange]     table[x expr=\thisrow{lvl}+0,y=average] {data/naca_ranks_67M_tol1e-3_PLU.dat};
        
                \addplot[mark=*,dashed,black]       table[x expr=\thisrow{lvl}+8,y=max] {data/naca_ranks_250k_tol1e-3_PLU.dat};
                \addplot[mark=*,dashed,red]         table[x expr=\thisrow{lvl}+6,y=max] {data/naca_ranks_1M_tol1e-3_PLU.dat};
                \addplot[mark=*,dashed,blue]        table[x expr=\thisrow{lvl}+4,y=max] {data/naca_ranks_4M_tol1e-3_PLU.dat};
                \addplot[mark=*,dashed,hcgreen]     table[x expr=\thisrow{lvl}+2,y=max] {data/naca_ranks_16M_tol1e-3_PLU.dat};
                \addplot[mark=*,dashed,hcorange]    table[x expr=\thisrow{lvl}+0,y=max] {data/naca_ranks_67M_tol1e-3_PLU.dat};
                \legend{Test 1,Test 2,Test 3,Test 4,Test 5};
            \end{semilogyaxis}
        \end{tikzpicture}
    }\\
    \subfloat[GMRES iterations]{
        \begin{tikzpicture}\label{gmresiteration}
            \begin{axis}[
            axis lines = left,
            ymin = 0,
                width=0.45\textwidth,xlabel={Nonlinear iteration},ylabel={GMRES iterations},
            ] 
                \addplot[black,mark=otimes]         table[x=step,y=gmres_steps,discard if not={N}{262144}]   {data/naca_plu_all_tol1e-3_gmres1e-3.dat};
                \addplot[red,mark=*]                table[x=step,y=gmres_steps,discard if not={N}{1048576}]  {data/naca_plu_all_tol1e-3_gmres1e-3.dat};
                \addplot[blue,mark=diamond]         table[x=step,y=gmres_steps,discard if not={N}{4194304}]  {data/naca_plu_all_tol1e-3_gmres1e-3.dat};
                \addplot[hcgreen,mark=square]       table[x=step,y=gmres_steps,discard if not={N}{16777216}] {data/naca_plu_all_tol1e-3_gmres1e-3.dat};
                \addplot[hcorange,mark=triangle]    table[x=step,y=gmres_steps,discard if not={N}{67108864}] {data/naca_plu_all_tol1e-3_gmres1e-3.dat};
            \end{axis}
        \end{tikzpicture}
    }
    \subfloat[Solve time]{
        \begin{tikzpicture}\label{solvetime}
            \begin{axis}[
            axis lines = left,
                ymin=0,width=0.45\textwidth,xlabel={Nonlinear iteration},ylabel={GMRES time{, s}},legend pos=outer north east,
            ] 
                \addplot[black,mark=otimes]         table[x=step,y=tgmres,discard if not={N}{262144}]   {data/naca_plu_all_tol1e-3_gmres1e-3.dat};
                \addplot[red,mark=*]                table[x=step,y=tgmres,discard if not={N}{1048576}]  {data/naca_plu_all_tol1e-3_gmres1e-3.dat};
                \addplot[blue,mark=diamond]         table[x=step,y=tgmres,discard if not={N}{4194304}]  {data/naca_plu_all_tol1e-3_gmres1e-3.dat};
                \addplot[hcgreen,mark=square]       table[x=step,y=tgmres,discard if not={N}{16777216}] {data/naca_plu_all_tol1e-3_gmres1e-3.dat};
                \addplot[hcorange,mark=triangle]    table[x=step,y=tgmres,discard if not={N}{67108864}] {data/naca_plu_all_tol1e-3_gmres1e-3.dat};
                \legend{Test 1,Test 2,Test 3,Test 4,Test 5};
            \end{axis}
        \end{tikzpicture}
    }
    \caption[Naca results]{Weak scaling results from 8 cores for Test 1 to 2048 cores for Test 5.}
    \label{fig-pspand:naca_plu_vs_pluq}
\end{figure}
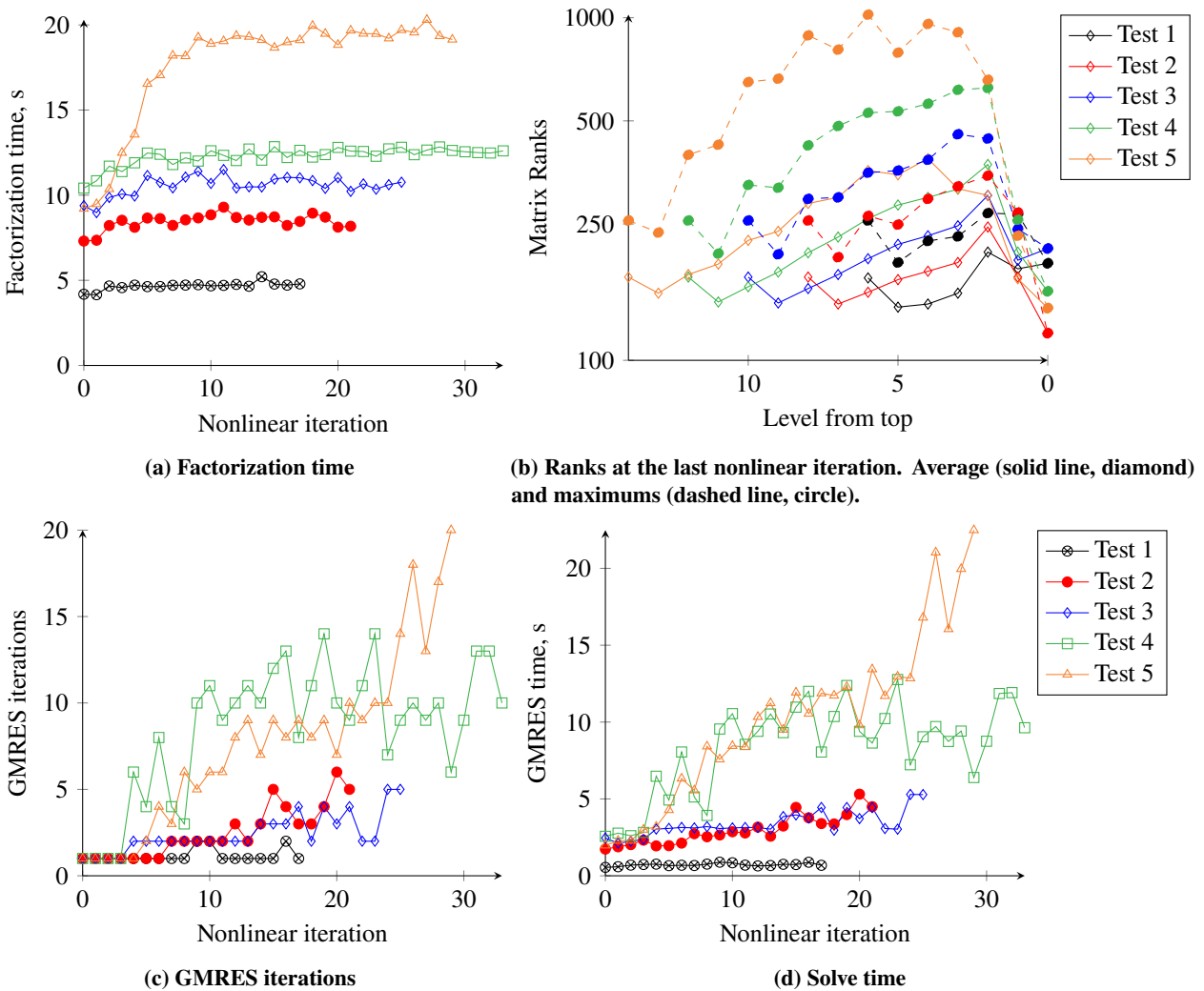

Fig. \ref{fig-pspand:naca_plu_vs_pluq} presents all the results throughout all nonlinear iterations. Here, we make a few remarks on the scalability of spaND. The factorization time of computing the preconditioner, spaND, scales well with the problem size, albeit with an increase in the larger test case (Fig. \ref{fact}). The average matrix ranks grow slowly with the problem size whereas the maximum matrix ranks grow at a higher rate and are not localized to the low-level leaves (Fig. \ref{ranks}). This can slow down spaND significantly since high ranks lead to longer sparsification time and high ranks towards the top of the hierarchical tree lead to poorer concurrency. We conjecture that the relatively high maximum ranks in larger test cases are degrading performance and causing the increase in factorization time. The number of GMRES steps remains very low at all problem sizes and does not vary much during the nonlinear convergence (Fig. \ref{gmresiteration}). The solve time also scales with the number of GMRES steps, indicating that applying the preconditioner scales well with the problem size (Fig. \ref{solvetime}).

\subsection{Scalability study of CAGMRES}\label{sec:cagmres}
As seen in the previous section, linear systems preconditioned by spaND require very few GMRES steps to reach convergence. The communication cost is therefore minimal, which justifies the use of a standard GMRES implementation. However, one could also adjust the tunable tolerance $\varepsilon$ in spaND to achieve better performance at the expense of more GMRES iterations. On large machines, the communication overhead becomes significant and it requires CA techniques to achieve scalability. 

For this purpose, we conduct an additional set of numerical experiments designed to stress the Krylov linear solver in order to demonstrate the scalability of CAGMRES. In this section, we focus on the scalability of the CAGMRES linear solver instead of the overall non-linear convergence. The tolerance parameter in spaND is increased to $10^{-1}$ and the tolerance of the CAGMRES solver is lowered to $10^{-10}$. Maximum number of linear iterations is increased to $100$ as well to collect meaningful performance data. Scalability tests are again in the form of weak scaling tests where the number of processing cores scale linearly with the problem size\footnote{Tests performed on a cluster equipped with dual-sockets and 24 cores Intel(R) Xeon(R) CPU E5-2680v3 @2.50GHz with 128GB of RAM per node.}. The statistics and tunable parameters of the tests used in this section are summarized in Table \ref{table:2}.

\begin{table}[!ht]
\caption{Summary of statistics and tunable paramters of additional tests used in Sec. \ref{sec:cagmres}}
\centering
\begin{tabular}{M{1cm}M{1.5cm}M{2cm}M{1.5cm}M{1.8cm}M{1cm}M{1.5cm}M{1.5cm}}
 \hline
 Test & No. of elements & Jacobian dimension, $N$ & No. of cores & SpaND tolerance, $\varepsilon$  & $\ell_{\max}$ & CAGMRES tolerance & CAGMRES max. steps \\
 \hline\hline
 6 & 65,536 & 1,048,576 & 48 & $10^{-1}$ & 12 & $10^{-10}$ & 100\\
 7 & 262,144 & 4,194,304 & 192 & $10^{-1}$ & 14 & $10^{-10}$ & 100\\
 8 & 1,048,576 & 16,777,216 & 768 & $10^{-1}$ & 16 & $10^{-10}$ & 100\\
 \hline
\end{tabular}
\label{table:2}
\end{table}

\begin{figure}[h]
    \centering
    \begin{tikzpicture}
    \begin{axis}
    [
        width=3.5in,
        height=3in,
        axis lines = left,
        xlabel = Linear iteration $j$,
        ylabel = $||Q_j^TQ_j-I_j||_F$,
        legend pos=north east,
        xtick={0, 4, 9, 14, 19, 24},
        xticklabels={4, 20, 40, 60, 80, 100},
        ymax=-9,
        ytick = {-10, -11, -12, -13, -14, -15},
        yticklabels={$10^{-10}$, $10^{-11}$, $10^{-12}$, $10^{-13}$, $10^{-14}$, $10^{-15}$},
    ]
    \addplot [black,mark=otimes]
    table{data/ortho_error_CGS1shift_linear.dat};
    \addplot [red,mark=*]
    table{data/ortho_error_CGS1noshift_linear.dat};
    \addplot [blue,mark=diamond]
    table{data/ortho_error_CGS2noshift_linear.dat};
    \addplot [hcgreen,mark=square]
    table{data/ortho_error_CGS2shift_linear.dat};
    \legend{$s$-step,$s$-step + Newton basis,$s$-step + Re-ortho., $s$-step + Newton basis + Re-ortho.};
    \end{axis}
    \end{tikzpicture}
    \caption{Orthogonality error $(||Q_j^TQ_j-I_j||_F)$ of CAGMRES with different stability measures for Jacobian matrix generated in Test 6 at nonlinear iteration 5. A step size of $s=4$ is used such that the linear iteration $j$ increases in increment of $s$.}
    \label{fig:stability}
\end{figure}
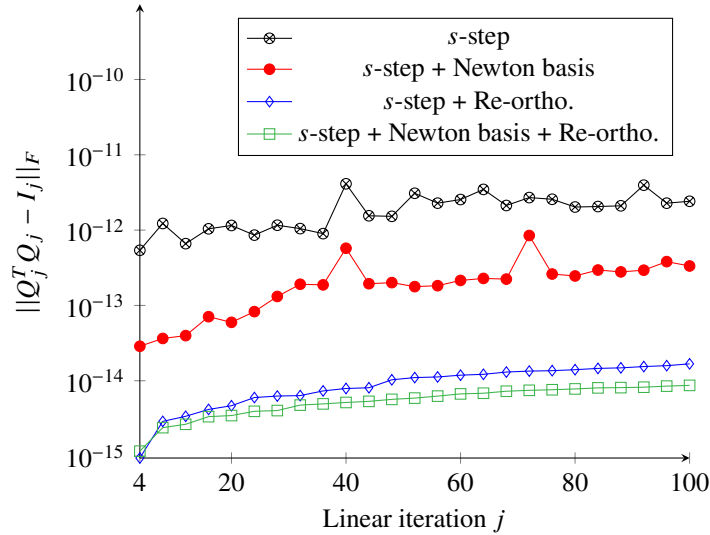

At each nonlinear iteration, the preconditioned linear system is solved using $s$-step CAGMRES with single-reduce re-orthogonalization scheme presented in Algorithm \ref{alg2}. A step size $s = 4$ is used to reduce the latency cost as only one global reduction is required for every $s$ basis vectors. The effective linear iteration count is reduced by a factor of $s$ correspondingly. Single-reduce re-orthogonalization scheme using two CGS steps and two CholQR steps is incorporated to tackle the numerical instability introduced by finite-precision implementation of $s$-step CAGMRES. 
This introduces additional computational cost but is found necessary as traditional block CGS in the original $s$-step CAGMRES implementation \cite{hoemmen2010communication} has an orthogonality error that is proportional to the square of the condition number of the $s$ block \cite{swirydowicz2020low,yamazaki2020low}. A Newton basis is included for the same stability concern \cite{bai1994newton}. The shifts used for the Newton basis are computed using the Ritz value from $s$ iterations of standard GMRES and are arranged in a Leja ordering for real arithmetic \cite{reichel1990newton,bai1994newton}. The effect of different stability measures on orthogonality error are illustrated in Fig. \ref{fig:stability}.


\begin{figure}[h]
    \centering
\begin{tikzpicture}
  \begin{semilogyaxis}[
    ybar stacked,
    bar shift=-8pt,
    legend pos = outer north east, 
    axis lines = left,
    xmin = -0.3,
    xmax = 2.3,
    ylabel = Wall time{, s},
    ymax=110,
    xtick = {0, 1, 2}, 
    xticklabels={Test 6, Test 7, Test 8},
    ytick = {10, 50, 100},
    yticklabels = {10, 50, 100},
    xticklabel style={align=center},
    ]
    \addplot [fill=red] table [x index = 0, y index = 1] {data/cagmres.dat};
    \addplot [fill=hcgreen] table [x index = 0, y index = 2] {data/cagmres.dat};
    \addplot [fill=black!50] table [x index = 0, y index = 3] {data/cagmres.dat};
    \addplot [fill=black] table [x index = 0, y index = 4] {data/cagmres.dat};
    \addplot [fill=blue] table [x index = 0, y index = 6] {data/cagmres.dat};
    \addplot [fill=hcorange] table [x index = 0, y index = 5] {data/cagmres.dat};
    \draw[thick, densely dashed] (-1,9.21) -- (3,9.21);
    \legend{SpMV, CGS/MGS, Re-ortho, Reduction, Solve, Other}
  \end{semilogyaxis}
  \begin{semilogyaxis}[
    ybar stacked,
    bar shift = 8pt,
    axis lines = left,
    xmin = -0.3,
    xmax = 2.3,
    ymax=110,
    xtick = {0, 1, 2}, 
    xticklabels={Test 6, Test 7, Test 8},
    ytick = {10, 50, 100},
    yticklabels = {10, 50, 100},
    xticklabel style={align=center},
    ]
    \addplot [fill = red] table [x index = 0, y index = 1] {data/gmres.dat};
    \addplot [fill = hcgreen] table [x index = 0, y index = 2] {data/gmres.dat}; 
    \addplot [fill = hcorange] table [x index = 0, y index = 3] {data/gmres.dat};
    \addplot [fill = black] table [x index = 0, y index = 4] {data/gmres.dat};
    \addplot [fill = blue] table [x index = 0, y index = 5] {data/gmres.dat};
  \end{semilogyaxis}
\end{tikzpicture}
\caption{Timing of various components of CAGMRES (left) and GMRES (right) for Tests 6-8 at nonlinear iteration 5. Step-size $s = 4$. (CGS/MGS: timing for arithmetic computation of CGS or MGS. Reduction: timing for MPI-based global reduction. Solve: timing for applying the preconditioner, spaND.)}
    \label{fig:cagmres}
\end{figure}
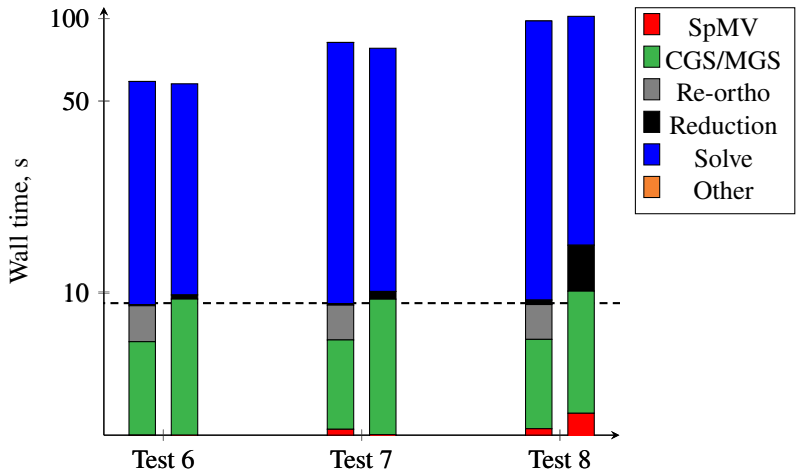

Using $s$-step CAGMRES with single-reduce re-orthogonalization schemes and Newton basis, the total computational time for 100 linear iterations of CAGMRES at each nonlinear iteration in each test is collected. Fig. \ref{fig:cagmres} shows the weak scaling results along with standard GMRES timings as baseline for comparison. We first highlight the near-constant timings of all components of CAGMRES (dashed line) except the solve step which corresponds to applying the preconditioner, spaND. Comparing to its GMRES counterpart, CAGMRES achieves scalability mainly due to its reduced number of global reductions which become more and more expensive on larger computing architecture. The orthogonalization step CGS takes less computational time than MGS as CGS can leverage on BLAS-3 matrix-matrix primitives that are more optimized than BLAS-2 matrix-vector kernels in MGS. The re-orthogonalization scheme introduces additional computational cost for stability but does not introduce extra communication cost. Applying the preconditioner takes about 80\% of total time in all cases but grow slowly with respect to the problem size.

\section{Conclusions}
\label{sec:conclusion}
This paper has described a new methodology for preconditioned scalable solvers for CFD workflows. The methodology consists of a hierarchical preconditioning strategy combined with communication-avoiding GMRES to enable scaling to very large numbers of processors in a weak-scaling sense. Both the preconditioning technique, spaND (with TaskTorrent), and the CAGMRES linear solver demonstrate near-linear weak scaling up to 2,048 cores in the context of a high-order DG solver within the SU2 framework. spaND shows scalable computational cost in the factorization of the Jacobian matrix while approximating the inverse with high accuracy, leading to constantly-small number of subsequent GMRES iterations. The CAGMRES solver developed minimizes communication overhead by employing the $s$-step technique to reduce latency and single-reduce re-orthogonalization scheme for stability. 

The set of numerical experiments presented in this work focuses on $h$-refinement while keeping the Jacobian block structure constant. Future efforts will extend the scope of the scalability analysis to include $p$-refinement as well. The final objective aims at demonstrating scalable, asymptotic performance for large-scale turbulent flows with complex geometries.

\section*{Acknowledgements}
\label{sec:acknowledgement}
Part of the computing resources for this runs executed in this paper were performed on the Stanford Research Computing Center cluster. We would like to thank Stanford University and the Stanford Research Computing Center for providing computational resources and support that contributed to these research results.
Léopold Cambier was supported by a fellowship from Total SE. Zan Xu was supported by NASA Grant 80NSSC18M0152 from the NASA Transformational Tools and Technologies program.

\bibliography{bibliography}

\end{document}